\documentclass[letterpaper, 10pt]{ieeeconf}     % Comment this line out

\usepackage{amssymb,amsmath,amsfonts,mathrsfs}
\usepackage{graphicx}

\newtheorem{theorem}{Theorem}[section]
\newtheorem{remark}{Remark}[section]

\newtheorem{definition}{Definition}[section]

\newcommand{\cH}{\mathcal{H}}

\newcommand{\cD}{\mathcal{D}}
\newcommand{\cC}{\mathcal{C}}
\newcommand{\cM}{\mathcal{M}}
\newcommand{\cV}{\mathcal{V}}

\newcommand{\bbR}{\mathbb{R}}

\newcommand{\bbE}{\mathbb{E}}

\newcommand{\bbP}{\mathbb{P}}

\renewcommand{\Re}[1]{\mathrm{Re}\{#1\}}

\newcommand{\iprod}[1]{\left\langle #1 \right\rangle}

\newcommand{\ket}[1]{\left|#1\right\rangle}
\newcommand{\kebr}[1]{\left|#1\right\rangle\left\langle #1 \right|}
\newcommand{\wto}{\hookrightarrow}
\DeclareMathOperator{\supp}{\mathrm{supp}}
\newcommand{\tr}[1]{\mathrm{Tr}\left\{#1 \right\}}

\DeclareMathOperator*{\argmin}{argmin}

\title{Approximate stabilization of an infinite dimensional quantum stochastic system}
\author{Ram Somaraju, Mazyar Mirrahimi\thanks{Ram Somaraju and Mazyar Mirrahimi are with INRIA Rocquencourt, Domaine
de Voluceau, B.P. 105, 78153 Le Chesnay cedex, France,
(ram.somaraju, mazyar.mirrahimi)@inria.fr}\thanks{Ram Somaraju and Mazyar Mirrahimi acknowledge support from ``Agence Nationale de la Recherche'' (ANR), Projet Jeunes Chercheurs EPOQ2 number ANR-09-JCJC-0070.} and Pierre Rouchon\thanks{P. Rouchon is with Mines ParisTech, Centre Automatique et Syst\'{e}mes, Math\'{e}matiques et Syst\'{e}mes, 60 Bd Saint Michel, 75272 Paris cedex 06,
France, pierre.rouchon@mines-paristech.fr}\thanks{Pierre Rouchon acknowledges support from ANR (CQUID).}}

\begin{document}
\maketitle


\begin{abstract}
We propose a feedback scheme for preparation of photon number states in a microwave cavity. Quantum Non-Demolition (QND) measurements of the cavity field and a control signal consisting of a microwave pulse injected into the cavity are used to drive the system towards a desired target photon number state. Unlike previous work, we do not use the Galerkin approximation of truncating the infinite-dimensional system Hilbert space into a finite-dimensional subspace. We use an (unbounded) strict Lyapunov function and prove that a feedback scheme that minimizes the expectation value of the Lyapunov function at each time step stabilizes the system at the desired photon number state with (a pre-specified) arbitrarily high probability. Simulations of this scheme demonstrate that we improve the performance of the controller by reducing ``leakage'' to high photon numbers. 
\end{abstract}
\section{Introduction}
Quantum Non-Demolition (QND) measurements have been used to detect and/or produce highly non-classical states of light in trapped super-conducting cavities~\cite{Deleglise2008,Gleyzes2007,Guerlin2007} (see~\cite[Ch. 5]{Haroche2006} for a description
of such quantum electro-dynamical  systems and~\cite{Brune1992} for detailed physical models with QND measures of light using atoms). In this paper we examine the feedback stabilization of such experimental setups near a pre-specified target photon number state. Such photon number states, with a precisely defined number of photons, are highly non-classical and have potential applications in quantum information and computation.

The state of the cavity may be described on a Fock space $\cH$, which is a particular type of Hilbert space that is used to describe the dynamics of a quantum harmonic oscillator (see e.g.~\cite[Sec 3.1]{Haroche2006}). The cannonical orthonormal basis for this Hilbert space consists of the set of Fock states $\{\ket{0},\ket{1},\ket{2},\ldots\}$. Physically, the state $\ket{n}$ corresponds to a cavity state with precisely $n$ photons. In this paper we study the possibility of driving the state of the system to some pre-specified target state $\ket{\bar{n}}$. The feedback scheme uses the so called measurement back action and a control signal, which is a coherent light pulse injected into the cavity, to stabilize the system at the target state with high probability.

Such feedback schemes for this experimental setup were examined previously in~\cite{Mirrahimi2010,Dotsenko2009}. The overall control structure used in~\cite{Mirrahimi2010} is a quantum adaptation of the observer/controller structure widely used for classical systems (see, e.g.~\cite[Ch. 4]{Kailath1980}). The observer part consists of a discrete-time quantum filter, based on the observed detector clicks, to estimate the quantum-state of the cavity field. This estimated state is then used in a state-feedback based on Lyapunov design, the controller part.

As the Hilbert space $\cH$ is infinite dimensional it is difficult to design feedback controllers to drive the system towards a target state  (because closed and bounded subsets of $\cH$ are not compact). In~\cite{Mirrahimi2010}, the controller was designed by approximating the underlying Hilbert space $\cH$ with a finite-dimensional Galerkin approximation $\cH_{N_{max}}$. Here, $\cH_{N_{max}}$ is the linear subspace of $\cH$ spanned by the basis vectors $\ket{0},\ket{1},\ldots,\ket{N_{max}}$ and $N_{max} \gg \bar{n}$, our target sate. Physically this assumption leads to an artificial bound $N_{max}$ on the maximum number of photons that may be inside the cavity. In this paper we wish to design a controller for the full Hilbert space $\cH$ without using the finite dimensional approximation. The need to consider the full Hilbert space is motivated by simulations (see Section~\ref{sec:simul}) which indicate that  using the controller designed on a finite dimensional approximation results in ``leakage'' to higher photon numbers with some finite probability.

Controlling infinite dimensional quantum systems have previously been examined in the deterministic setting without measurements. Various approaches have been used to overcome the non-compactness of closed and bounded sets. One approach consists of proving approximate convergence results which show convergence to a neighborhood of the target state~\cite{Beauchard2009,Mirrahimi2009}. Alternatively, one examines weak convergence for example, in~\cite{Beauchard2010}. Other approaches such as using strict Lyapunov functions or strong convergence under restrictions on possible trajectories to compact sets have also been used in the context of infinite dimensional state-space for example in~\cite{Coron1998,Coron2007}.

The situation in our paper is different in the sense that the system under consideration is inherently stochastic due to quantum measurements. The system may be described using a discrete time Markov process on the set of unit vectors in the system Hilbert space as explained in Section~\ref{sec:sysDisc}. We use a strict Lyapunov function that restricts the system trajectories with high probability to compact sets as explained in Section~\ref{sec:mainRes}. We use the properties of weak-convergence of measures to show approximate convergence (i.e. with probability of convergence approaching one) of the discrete time Markov process towards the target state. 

We use a similar overall feedback scheme that is used in~\cite{Mirrahimi2010}. The entire feedback system is split into an observer part, a quantum filter, and a controller part based on a Lyapunov function. The quantum filter used to estimate the state is identical to the one used in~\cite{Mirrahimi2010} and we do not discuss the filter further in this paper. However we do not use the Galerkin approximation to design the controller. We show in Theorem~\ref{the:mainRes} that given any $\epsilon > 0$, we can drive our system to the target state $\bar{n}$ with probability greater than $1-\epsilon$. Simulations (see Section~\ref{sec:simul}) indicate that this controller provides improved performance with lower probability of having trajectory escaping towards infinite photon numbers. The precise choice of Lyapunov function is motivated by~\cite{Amini2011} that uses a similar form of the Lyapunov function in a finite dimensional setting. 
\subsection{Outline}
The remainder of the paper is organised as follows: in the following Section we describe the experimental setup and the Markovian jump dynamics of the system state. In Section~\ref{sec:mainRes} we state the main result of our paper including an outline of the proof of Theorem~\ref{the:mainRes}. We then present our simulation results in Section~\ref{sec:simul} and then our conclusions in the final Section.
\section{System description}\label{sec:sysDisc}
\begin{figure}
 \includegraphics[width=0.48\textwidth]{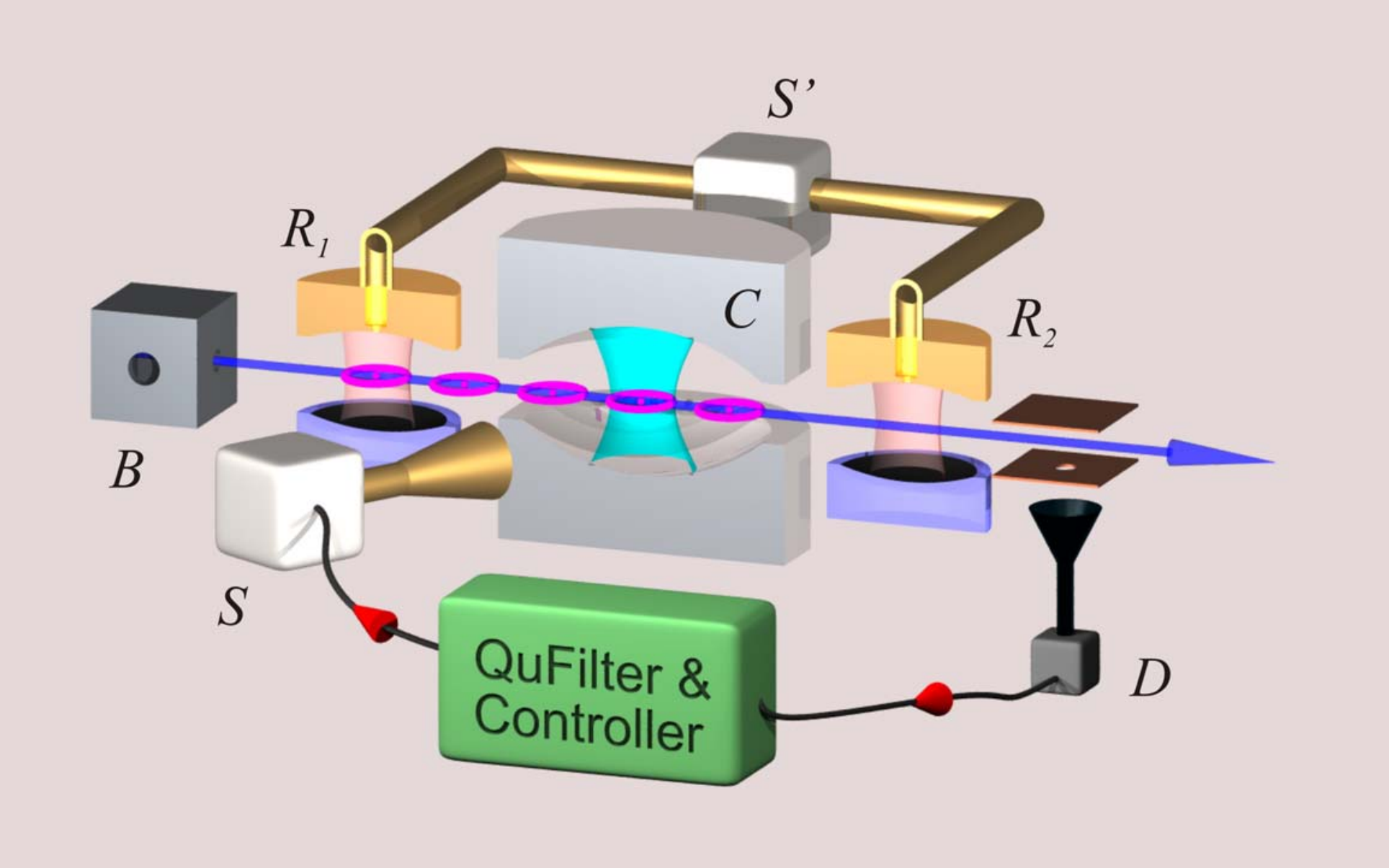}\\
\caption{The microwave cavity QED setup with its feedback scheme (in
green).}\label{fig:sys}
\end{figure}
The system, illustrated in Figure~\ref{fig:sys}, consists of 1) a high-$Q$ microwave cavity $C$, 2) an atom source $B$ that produces Rydberg atoms, 3) two low-Q Ramsey cavities $R_1$ and $R_2$, 4) an atom detector $D$ and 5) a microwave source $S$. The system may be modeled by a discrete-time Markov process, which takes into account the backaction of the measurement process (see e.g.~\cite[Ch. 4]{Haroche2006} and~\cite{Mirrahimi2010}).

Rydberg Atoms are sent from $B$, interact with the cavity $C$, entangling the state of the atom with that of the cavity and are then  detected in $D$. Each time-step, indexed by the integer $k$, corresponds to atom number $k$ crossing the cavity and interacting with the cavity. The state of the cavity in time step $k$ is described by a unit vector $\ket{\psi_k}\in \bar{B}_1$ for $k = 1,2,\ldots$. Here, $\bar{B}_1 = \{\ket{\psi}\in \cH:\|\ket{\psi}\| = 1\}$ is the set of possible cavity states. The change of the cavity state $\ket{\psi_k}$ at time-step $k$ to the state $\ket{\psi_{k+1}}$ at time-step $k + 1$ consists of two parts corresponding to the projective measurement of the cavity state, by detecting the state of the Rydberg atom in detector $D$ and also due to an appropriate coherent pulse (the control) injected into C.

Let $a$ and $a^\dagger$ be the photon annihilation and creation operators where $a\ket{n} = \sqrt{n} \ket{n-1}$ and $a^\dagger$ is the Hermition conjugate of $a$. Also, let $N = a^\dagger a$ be the diagonal number operator satisfying $N\ket{n} = n\ket{n}$. Let $D_\alpha = \exp(\alpha(a^\dagger - a))$ be the displacement operator which is a unitary operator that corresponds to the  input of a coherent control field of amplitude $\alpha$ that is injected into the cavity. The amplitude $\alpha$ of the coherent field is the control that is used to manipulate the system. Let $M_g = \cos(\theta + N\phi)$ and $M_e = \sin(\theta + N\phi)$ be the measurement operators, where $\theta$ and $\phi$ are experimental parameters. Physically, the measurement operator $M_s$, $s\in \{e,g\}$ correspond to the state of the detected atom in either the ground state $\ket{g}$ or the excited state $\ket{e}$ .
%So, if the atom $k$ is detected in the ground state $\ket{g}$ and then the state of the cavity at time step $k+1$ is $\frac{M_g\ket{\psi}}{\|M_g\ket{\psi}\|}$ given state $\ket{\psi}$ at time step $k$.

We model these dynamics by a Markov process
\begin{eqnarray}
\ket{\psi_{k+1/2}} &=&\frac{ M_s \ket{\psi_k}}{\|M_s\ket{\psi_k}\|} \textrm{ with prob. } \|M_s\ket{\psi_k}\|^2\label{eqn:half}\\
\ket{\psi_{k+1}} &=& D_{\alpha_k} \ket{\psi_{k+1/2}}.\label{eqn:one}
\end{eqnarray}
Here $s\in \{e,g\}$ and the control $\alpha_k\in \bbR$.
\begin{remark}
The time evolution from the step $k$ to $k + 1$, consists of two types of evolutions: a projective measurement by the operators  $M_s$ and a coherent injection involving operator $D_\alpha$. For the sake of simplicity, we will use the notation of $\ket{\psi_{k+1/2}}$ to illustrate this intermediate step.
\end{remark}
\begin{remark}
Let $\cM_1$ be the set of all probability measures on $\bar{B}_1$. Then the Equations~\eqref{eqn:half} and~\eqref{eqn:one} determine a stochastic flow in $\cM_1$ and we denote by $\Gamma_k(\mu_0)$ the probability distribution of $\ket{\psi_k}$, given $\mu_0$, the probability distribution of $\ket{\psi_0}$.
\end{remark}
\section{Global (approximate) feedback stabilization}\label{sec:mainRes}
We wish to use the control $\alpha_k$ to drive the system into a pre-specified target state $\ket{\bar{n}}$ with high probability. That is, we wish to show that the sequence $\Gamma_k(\mu)$ converges to the set of probability measures $\Omega_\infty$ where for all $\mu_\infty\in \Omega_\infty$, $\mu_\infty(\ket{n})$ is big.

In order to achieve this we use a Lyapunov function~\eqref{eqn:lyap} and at each time step $k$ we choose the feedback control $\alpha_k$ to minimize the Lyapunov function. Before discussing the choice of the Lyapunov function in Subsection~\ref{subsec:lyap} we recall some facts concerning the convergence of probability measures
\subsection{Convergence of probability measures}
We refer the interested reader to~\cite{Merkle2000,Billingsley1999} for results pertaining to convergence of probability measures. We denote by $\cC$ the set of all continuous bounded functions on $\bar{B}_1$.
\begin{definition}\label{defn:convMeas} We say that a sequence of probability measure $\{\mu_n\}_{n=1}^\infty\subset \cM_1$ converges (weak-$\ast$) to a probability measure $\mu\in \cM_1$ if for all $f\in \cC$
\begin{equation*}
\lim_{n\to\infty} \bbE_{\mu_n}[f] = \bbE_\mu[f]
\end{equation*}
and we write
\begin{equation*}
\mu_n\wto \mu.
\end{equation*}
\end{definition}
It can be shown that if $\mu_n\wto\mu_\infty$ then for all open sets $W$, 
\begin{equation}\label{eqn:open}
\liminf_{n\to \infty} \mu_n(W) \geq \mu_\infty(W).
\end{equation}
A set of probability measures $S\subset \cM_1$ is said to be \emph{tight}~\cite[p. 9]{Billingsley1999} if for all $\epsilon > 0$ there exists a compact set $K_\epsilon \subset \bar{B}_1$ such that for all $\mu\in S$, $\mu(K_\epsilon) > 1- \epsilon$.
\begin{theorem}[Prohorov's theorem]\label{the:prohorov}
Any tight sequence of probability measures has a (weak-$\ast$) converging subsequence.
\end{theorem}
We also recall Doob's inequality. Let $X_n$ be a Markov process on some state space X. Suppose that there is a non-negative function $V(x)$ satisfying $\bbE[V(X_1)|X_0 = x)]- V (x) \leq 0$, then Doob's inequality states
\begin{equation}\label{eqn:doob}
\bbP\left(\sup_{n\geq 0}V(X_n) \geq \gamma|X_0 = x\right) \leq \frac{V(x)}{\gamma}.
\end{equation}
\subsection{Lyapunov function and control signal $\alpha_k$}\label{subsec:lyap}
We now introduce our Lyapunov function $V$ and explain the intuition behind this peculiar form of this function. The function, $V:\bar{B}_1\to [0,\infty]$ is defined
\begin{eqnarray}
V(\ket{\psi}) &=& \sum_{n=0}^\infty \sigma_n\left|\iprod{\psi|n}\right|^2 + \delta(\cos^4(\phi_{\bar{n}}) + \sin^4(\phi_{\bar{n}})) \nonumber\\
&& -\: \delta\big( \left\|M_g\ket{\psi} \right\|^4 + \left\|M_e \ket{\psi}\right\|^4\big).\label{eqn:lyap}
\end{eqnarray}
Here
\begin{equation*}
\phi_{n} = \theta + n\phi ,
\end{equation*}
$\delta > 0$ is a small positive number and
\begin{equation}\label{eqn:sigmaDef}
\sigma_n = \left\{\begin{array}{ll}
\frac{1}{8} +\sum_{k=1}^{\bar{n}} \frac{1}{k} - \frac{1}{k^2},& \textrm{ if $n=0$}\\
\sum_{k=n+1}^{\bar{n}} \frac{1}{k} - \frac{1}{k^2},& \textrm{ if $1 \leq n < \bar{n}$}\\
0, &  \textrm{ if $n = \bar{n}$ }\\
\sum_{k=\bar{n}+1}^n \frac{1}{k} + \frac{1}{k^2},& \textrm{ if $n > \bar{n}$}
\end{array}\right.
\end{equation}
We set $\cD(V)\subset \bar{B}_1$ to be the set of all $\ket{\psi}\in \bar{B}_1$ where the above Lyapunov function  is finite. We note that coherent states, which are states that are of relevance in practical experiments are in $\cD(V)$.

We choose a feedback that minimizes the expectation value of the Lyapunov function in every time-step $k$. Indeed, applying the result of the $k$'th measurement, we know the state $\ket{\psi_{k+1/2}}$ and we choose $\alpha_k$ as follows 
\begin{equation}\label{eqn:control}
\alpha_k = \argmin_{\alpha\in [-\bar{\alpha},\bar{\alpha}]}  V\left(D_{\alpha}\ket{\psi_{k+1/2}}\right)
\end{equation}
for some positive constant $\bar{\alpha}$.
\begin{remark} The Lyapunov function is chosen  to be this specific form to serve three purposes -
\begin{enumerate}
\item We choose the sequence $\sigma_n \to \infty$ as $n\to \infty$. This guarantees that if we choose $\alpha_k$ to minimize the  expectation value of the Lyapunov function then the trajectories of the Markov process are restricted to a compact set in $\bar{B}_1$ with probability arbitrarily close to 1. This implies that the $\omega$-limit set of the process is non-empty  (see Step 2 in the Proof of Theorem~\ref{the:mainRes}).
\item The term $-\delta(\|M_g\ket{\psi}\|^4 +\|M_g\ket{\psi}\|^4)$ is chosen such that the Lyapunov function is a strict Lyapunov functions for the Fock states. This implies that the support of the $\omega$-limit set only contains Fock states (see Step 3 in the Proof of Theorem~\ref{the:mainRes}).
\item The relative magnitudes of the coefficients $\sigma_n$ have been chosen such that $V(\ket{\bar{n}})$ is a strict global minimum of $V$. Moreover given any $M > \bar{n}$ we can choose $\delta,\bar{\alpha}$ such that for all $M \geq m\neq \bar{n}$, and for all $\ket{\psi}$ in a neighborhood of $\ket{m}$, $V(D_\alpha\ket{\psi})$ does not have a local minimum at $\alpha = 0$. This implies that if $\ket{\psi_k}$ is in this neighborhood of $\ket{m}$ then we can choose an $\alpha_k\in [-\bar{\alpha},\bar{\alpha}]$ to decrease the Lyapunov function and move $\ket{\psi_k}$ away from $\ket{m}$ by some finite distance with probability 1  (see Steps 4 and 5 in the Proof of Theorem~\ref{the:mainRes}).
\end{enumerate}
\end{remark}
\subsection{Main Result}
We make the following assumption.
\begin{itemize}
\item[A1] The eigenvalues of $M_g$ and $M_e$ are non-degenerate. This is equivalent to the assumption that $\pi/\phi$ is not a rational number.
\end{itemize}
The quantum filter uses the statistics of the measurement of whether the atom is in the ground or excited state to estimate the cavity's state. Therefore if one of the eigenvalues of $M_g$ (or $M_e$) is degenerate then the measurement statistics will be the same for more than one photon number state. Therefore it is not possible to control the system effectively in this case (However, as explained in Remark~\ref{rem:weakAss} below, we may weaken this assumption slightly).

The following Theorem is our main result.
\begin{theorem}\label{the:mainRes}
If we assume $A1$ to be true then given any $\epsilon > 0$ and $C > 0$, there exist constants $\delta >0$ and $\bar{\alpha}$ such that for all $\mu$ satisfying $\bbE_\mu[V] \leq C$, $\Gamma_n(\mu)$ converges to a limit set $\Omega$. Moreover for all $\mu_\infty\in \Omega$, $\ket{\psi} \in \supp(\mu_\infty)$ only if $\ket{\psi}$ is one of the Fock states $\ket{n}$ and
\begin{equation*}
\mu_\infty(\{\ket{\bar{n}}\}) \geq 1 - \epsilon.
\end{equation*}
\end{theorem}
The proof is split into 5 steps:
\begin{enumerate}
\item $V(\ket{\psi_k})$ is a super-martingale that is bounded from below.
\item The sequence of measures $\Gamma_k(\mu)$ is tight and therefore has a converging subsequence. Hence the set $\Omega$ is non-empty.
\item If $\Gamma_{k_l}(\mu) \to \mu_\infty$ then the support set of $\mu_\infty$ only consists of Fock states.
\item Let $M',C' > 0$ be given. Then for all $M' \geq m\neq \bar{n}$, $\delta$ and $\bar{\alpha}$ may be chosen small enough such that for $\kappa > 0$ small enough and all $\ket{\psi}$ in the neighborhood
\begin{equation}\label{eqn:vkm}
\cV^\kappa_m = \{\ket{\psi}:\|\ket{\psi}-\ket{m}\| < \kappa, V(\ket{\psi}) > V(\ket{m}) - \kappa\}
\end{equation}
of $\ket{m}$, satisfying $V(\ket{\psi}) < C'$, we have for $|\alpha| < \bar{\alpha}$ the polynomial approximation
\begin{equation*}
V(D_\alpha\ket{\psi}) =\sum_{i = 0}^2 \frac{\alpha^i}{i!}f_i(\ket{\psi}) + O(\bar{\alpha}^3) + O(\delta)
\end{equation*}
and $f_2(\ket{\psi}) < \gamma < 0$ for some constant $\gamma$. The term $O(\bar{\alpha}^3)$ only depends on $C'$ and not on $\ket{\psi}$ and the term $O(\delta)$ is independent of both $\ket{\psi}$ and $C'$.
\item Because $\gamma$ is negative, we can choose $\bar{\alpha}$ and $\delta$ small enough such that the probability of convergence to the Fock states $\ket{m}$ for $m \neq \bar{n}$ may be made arbitrarily small. Therefore
\begin{equation*}
\mu_\infty(\ket{\bar{n}}) = 1 - \sum_{\substack{m=0\\m\neq\bar{n}}}^\infty\mu_\infty(\ket{m})
\end{equation*}
may be made arbitrarily big.
\end{enumerate}

Below we sketch the proofs of each of the above steps. The interested reader is referred to~\cite{Somaraju2011} for further details on the proof which are beyond the scope of a short note.
\begin{proof}[Proof of step 1]
We can write
\begin{equation*}
\bbE\left[V(\ket{\psi_{k+1}})\big|\ket{\psi_k}\right] - V(\ket{\psi_k}) = K_1(\ket{\psi_k}) + K_2(\ket{\psi_k})
\end{equation*}
where,
\begin{eqnarray}
K_1(\ket{\psi_k}) &\triangleq& \min_{\alpha\in [-\bar{\alpha},\bar{\alpha}]} \bbE\left[V\left(D_{\alpha}\left(\ket{\psi_{k+1/2}}\right)\right)\big|\ket{\psi_k}\right] \nonumber\\
&& -\: \bbE\left[V(\ket{D_0(\psi_{k+1/2}})\big|\ket{\psi_k}\right],\label{eqn:k1}\nonumber\\
K_2(\ket{\psi_k}) &\triangleq&  \bbE\left[V\left(D_{0}\left(\ket{\psi_{k+1/2}}\right)\right)\big|\ket{\psi_k}\right] \nonumber\\&&- V(\ket{\psi_k}).\label{eqn:k2}
\end{eqnarray}
It is obvious that $K_1(\ket{\psi}) \leq 0$ and after simple but tedious manipulations, we get
\begin{equation}\label{eqn:k2}
K_2(\ket{\psi}) = \frac{-2\left(\|M_g^2\ket{\psi}\|^2 - \|M_g\ket{\psi}\|^4\right)^2}{\tr{M_g^2\rho}\tr{M_e^2\rho}} \leq 0.
\end{equation}
Therefore, $V(\psi_k)$ is a super-martingale.
\end{proof}
\begin{proof}[Proof of step 2]
Let $\epsilon> 0$ be given. Because $V(\ket{\psi_k})$ is a supermartingale, Doob's inequality~\eqref{eqn:doob} gives us
\begin{equation}\label{eqn:doob}
\bbP\left(\sup_{k\geq 0}V(\ket{\psi_k}) \geq \frac{\bbE_\mu[V]}{\epsilon}\right) \leq \epsilon.
\end{equation}
If we set,
\begin{equation*}
K_\epsilon = \{\ket{\psi}:V(\ket{\psi}) \leq \bbE_\mu[V]/\epsilon\})
\end{equation*}
then for all $k> 0$, $[\Gamma_k(\mu)](K_\epsilon) > 1 - \epsilon$. Because, the sequence $\sigma_n \to \infty$ as $n\to \infty$, the set $K_\epsilon$ can be shown to be pre-compact in $\cH$. We can now apply Prohorov's Theorem~\ref{the:prohorov} to show that $\Gamma_n(\mu)$ has a converging subsequence. Therefore the limit set $\Omega = \{\mu_\infty\in \cM_1:\Gamma_{k_l}(\mu)\wto \mu_{\infty}\}$ is non-empty.
\end{proof}
\begin{proof}[Proof of step 3] Suppose some subsequence of $\Gamma_k(\mu)$ converges to $\mu_\infty\in \Omega$. From  step 1 we have $K_1(\ket{\psi_k}) + K_2(\ket{\psi_k}) \to 0$ as $k\to \infty$ and  because $K_1$ and $K_2$ are both non-negative we have
\begin{equation*}
\lim_{k\to\infty} \bbE_{\Gamma_k(\mu)}[K_2] = 0.
\end{equation*}
But, from~\eqref{eqn:k2} and the boundedness of $M_g$ and $M_e$, we know that $K_2$ is a continuous function on $\cH$. Therefore from Definition~\ref{defn:convMeas} of (weak-$\ast$) convergence of measures we get
\begin{equation}\label{eqn:expK2}
\bbE_{\mu_\infty}[K_2] = 0.
\end{equation}
But $K_2(\ket{\psi}) = 0$ implies $\|M_g^2\ket{\psi}\|^2 = \|M_g\ket{\psi}\|^4$. The Cauchy-Schwartz inequality gives
\begin{eqnarray*}
\|M_g^2\ket{\psi}\|^2 &=& \|M_g^2\ket{\psi}\|^2 \|\ket{\psi}\|^2\\
&=& \iprod{\psi M_g^2|M_g^2\psi}\cdot\iprod{\psi|{\psi}} \\
&\geq& |\iprod{\psi|M_g^2{\psi}}|^2\\
&=& \|M_g\ket{\psi}\|^4.
\end{eqnarray*}
with equality if and only if $\ket{\psi}$ and $M^2_g\ket{\psi}$ are co-linear. Therefore $K_2(\ket{\psi}) = 0$ implies (by Assumption$A1$) that $\ket{\psi}$ is a Fock state. Hence from~\eqref{eqn:expK2} we can conclude that the support set of $\mu_\infty$ only consists of the set of Fock states.
\end{proof}
\begin{proof}[Proof of step 4]

Set
\begin{equation*}
\hat{V}(\ket{\psi}) \triangleq \sum_{n=0}^\infty\sigma_n|\iprod{D_\alpha \psi|m}|^2
\end{equation*}
It can be shown~\cite{Somaraju2011} that $\hat{V}(D_\alpha\ket{\psi})$ is an analytic function of $\alpha$ if $\ket{\psi}$ satisfies $\hat{V}(\ket{\psi}) < \infty$. Moreover, for all $\ket{\psi}$ satisfying $\hat{V}(\ket{\psi}) < C'$ we have the second order polynomial approximation
\begin{equation*}
\hat{V}(D_\alpha\ket{\psi}) =\sum_{i = 0}^2 \frac{\alpha^i}{i!} \nabla_\alpha^i \hat{V}(D_\alpha\ket{\psi})\big|_{\alpha = 0} + O(\bar{\alpha}^3)
\end{equation*}
for all $|\alpha| < \bar{\alpha}$. In particular the $O(\bar{\alpha})$ term only depends on $C'$ and is independent of $\ket{\psi}$. Here $\nabla_\alpha^i (\cdot)|_{\alpha = 0}$ is the $i^{th}$ derivative of $(\cdot)$ w.r.t. $\alpha$ evaluated at $\alpha = 0$.

If we let $\ket{\psi} = \sum_{n=0}^\infty c_n\ket{n}$ and recall that $D_\alpha = \exp(\alpha(a-a^\dagger))$ then after some manipulations, we get
\begin{eqnarray*}
\lefteqn{\nabla_\alpha^2\hat{V}(D_\alpha\ket{\psi})\big|_{\alpha = 0}} \\&& = \sum_{n=0}^\infty |c_n|^2\big((n+1)\sigma_{n+1} + n\sigma_{n-1} - (2n+1)\sigma_n\big)\nonumber\\
&&+\:\Re{c_{n-1}c_{n+1}^\ast}\sqrt{n(n+1)}(\sigma_{n-1} + \sigma_{n+1} - 2\sigma_n).
\end{eqnarray*}

If $n\neq \bar{n}$ and $n\geq 2$ we have
\begin{equation*}
(n+1)\sigma_{n+1} + n\sigma_{n-1} - (2n+1)\sigma_n = \frac{-1}{n(n+1)}
\end{equation*}
and for $n= 0,1$ we get
\begin{equation*}
(n+1)\sigma_{n+1} + n\sigma_{n-1} - (2n+1)\sigma_n = \frac{-1}{4}
\end{equation*}

For any Fock state $\ket{m}$ with $m\neq \bar{n}$, $c_n = \delta_{mn}$, where $\delta_{mn}$ is the Kronecker-delta function and we have
\begin{equation*}
\nabla_\alpha^2\hat{V}(D_\alpha\ket{m})\big|_{\alpha = 0} = -\frac{1}{m(m+1)} < 0.
\end{equation*}
Because the terms $\sum_n|c_n|^2$ and $\sum_n \Re{c_{n+1}c_{n-1}^\ast}$ are bounded by the $\|\cdot\|$-norm in $\cH$, it can be shown that for $\kappa$ small enough we have $\nabla_\alpha^2\hat{V}(D_\alpha\ket{\psi})\big|_{\alpha = 0} < -\frac{1}{2m(m+1)}$ in the neighborhood $\cV^\kappa_m$ of $\ket{m}$, where $\cV_m^\kappa$ is given as in Equation~\eqref{eqn:vkm}.

But,
\begin{equation*}
\nabla_\alpha^2 V(D_\alpha\ket{\psi})\big|_{\alpha = 0} = \nabla_\alpha^2\hat{V}(D_\alpha\ket{\psi})\big|_{\alpha = 0} + O(\delta).
\end{equation*}
Hence, given any $M > \bar{n}$, step 4 above is true with $\gamma = -\frac{1}{2M(M+1)}$.
\end{proof}
\begin{proof}[Proof of step 5]
Let $\epsilon > 0$ be given. We show that $\mu_\infty(\{\ket{\bar{n}}\}) \geq 1-\epsilon$. From step 3 we know that the support of $\mu_\infty$ only consists of Fock states. Therefore using~\eqref{eqn:open}, we only need to show that there exists an open neighborhood $W$ of $\{\ket{m}:m\neq \bar{n}\}$ such that for $k$ big enough the $[\Gamma_k(\mu)](W) \leq \epsilon$.

We construct the set $W$ using two disjoint parts $W_1$ and $W_2$. We first show that there exists a $M$ big enough and a neighborhood $W_1$ of $\{\ket{M},\ket{M+1},\ldots\}$ such that $[\Gamma_k(\mu)](W_1) \leq \epsilon/2$ for all $k$. We then construct a neighborhood $W_2$ of $\{\ket{m}:0\leq m < M, m\neq \bar{n}\}$ such that $[\Gamma_k(\mu)](W_2) < \epsilon/2$ for $k$ large enough.

\paragraph{Construction of $W_1$}
Because $\sigma_m\to \infty$ there exists an $M$ large enough such that for all $m > M$, $\sigma_m > \frac{C}{\epsilon/4}$. We can choose a small enough neighborhood $W_1$ of $\{\ket{M},\ket{M+1},\ldots\}$ such that for all $\ket{\psi}$ in this neighborhood,
\begin{equation*}
V(\ket{\psi}) \geq \frac{\sigma_M}{2} \geq \frac{C}{\epsilon/2}
\end{equation*}
Because $\bbE_\mu[V] \leq C$, Doob's inequality implies the probability of $V(\ket{\psi_k}) > C/(\epsilon/2)$ is less than $\epsilon/2$. Therefore,
\begin{equation}\label{eqn:e1}
[\Gamma_k(\mu)](W_1) \leq \frac{\epsilon}{2}.
\end{equation}
\paragraph{Construction of $W_2$} We show that for $\kappa$ small enough we can choose 
\begin{equation*}
W_2 = \bigcup_{\substack{m=0\\m\neq\bar{n}}}^{M-1} \cV^\kappa_m
\end{equation*}
where $\cV_m^\kappa$ is as in~\eqref{eqn:vkm}.

From Doob's inequality, we have
\begin{equation}\label{eqn:e2}
[\Gamma_k(\mu)]\left(\left\{\ket{\psi}:V(\ket{\psi}) > \frac{C}{\epsilon/2}\right\}\right)\leq \epsilon/2.
\end{equation}
for all $k$. Therefore we can complete the proof if we show that for $\kappa$ small enough
\begin{equation*}
\lim_{k\to\infty} [\Gamma_k(\mu)](\hat{\cV}^\kappa_m) = 0,
\end{equation*} 
where
\begin{equation*}
\hat{\cV}^\kappa_m = \cV^\kappa_m\cap \left\{\ket{\psi}: V(\ket{\psi}) \leq \frac{C}{\epsilon/2}\right\}.
\end{equation*}

In Step 4 we set $C' = C/(\epsilon/2)$ and $M' = M$ and let $\kappa$ be small enough so that $\cV_m^\kappa$ is as given in step 4. Then, because $\gamma < 0$, we can choose $\bar{\alpha}$ and $\delta$ small enough so that there exists a constant $c> 0$ such that for all $\ket{\psi}\in\hat{ \cV}^\kappa_m$, $V(D_\alpha\ket{\psi}) - V(\ket{\psi}) < -c$, for some $\alpha \in [- \bar{\alpha},\bar{\alpha}]$. Because $M_g$ and $M_e$ are bounded operators and $M_g\ket{m}= \ket{m}$ and $M_e\ket{m} = \ket{m}$, $\kappa'$ can be chosen small enough such that if $\ket{\psi_k}\in \hat{\cV}^{\kappa'}_m$ then $\ket{\psi_{k+1/2}} \in \hat{\cV}_m^\kappa$ with probability $1$.

We claim that $\Gamma_k(\mu)(\hat{\cV}_m^{\kappa'/2})\to 0$. To see this note that because for all $\ket{\psi_k}\in \hat{\cV}_m^{\kappa'}$, $\ket{\psi_{k+1/2}}\in \hat{\cV}^\kappa_m$ with probability 1 and this implies
\begin{equation*}
V(\ket{\psi_{k+1}}) - V(\ket{\psi_k}) < -c\textrm{ with probability $1$}.
\end{equation*}
So if $\ket{\psi_k}\in \hat{\cV}_m^{\kappa'/2}$ then the Markov process is outside the set $\cV_m^{\kappa'}$ within a finite number of steps less than $\lceil C'/c\rceil$ with probability $1$. So if $\mu_k(\hat{\cV}_m^{\kappa'/2})$ does not approach zero, then the Markov process must enter the set $\hat{\cV}_m^{\kappa'/2}$ from outside the set ${\cV}_m^{\kappa'}$ infinitely many times. But by Doob's inequality~\eqref{eqn:doob} the probability of this happening once is less than $(V(\ket{m})-\kappa')/(V(\ket{m})-\kappa'/2) < 1$. Therefore the probability of this happening infinitely many times is zero. Thus $\lim_{k\to\infty}[\Gamma_k(\mu)](\hat{\cV}_m^{\kappa'/2}) = 0$. This combined with~\eqref{eqn:e1} and~\eqref{eqn:e2} gives, $\mu_\infty\{\ket{m}:m\neq \bar{n}\}\leq \epsilon$.

Therefore,
\begin{equation*}
\mu_\infty(\ket{\bar{n}}) = 1 - \sum_{\substack{m=0\\m\neq\bar{n}}}^\infty\mu_\infty(\ket{m}) \geq 1-\epsilon.
\end{equation*}
\end{proof}

\begin{remark}\label{rem:weakAss}
In step 2 we show that the only vectors in the support of $\mu_\infty$ are those corresponding to eigenvector of $M_s$. We then used assumption $A1$ to claim that the only eigenvectors of $M_s$ are the Fock states. We can however weaken this assumption to the following: eigenvalues corresponding to eigenvectors $\ket{m}$, $m< M$ are non-degenerate. This is because, we can show that if some eigenvector $\ket{\psi}$ is in the span of the set $\{\ket{M},\ket{M+1},\ldots\}$ then using the same argument as that used for $\ket{m}, m> M$, we can show that the probability of $\ket{\psi}$ is small. This is significant for cases where $M_g$ is a more complicated non-linear function of $N$, as is the case in a practical system.
\end{remark}
\section{Simulations}\label{sec:simul}
\begin{figure}
  \includegraphics[trim = 0mm 85mm 0mm 90mm, clip, width=0.45\textwidth]{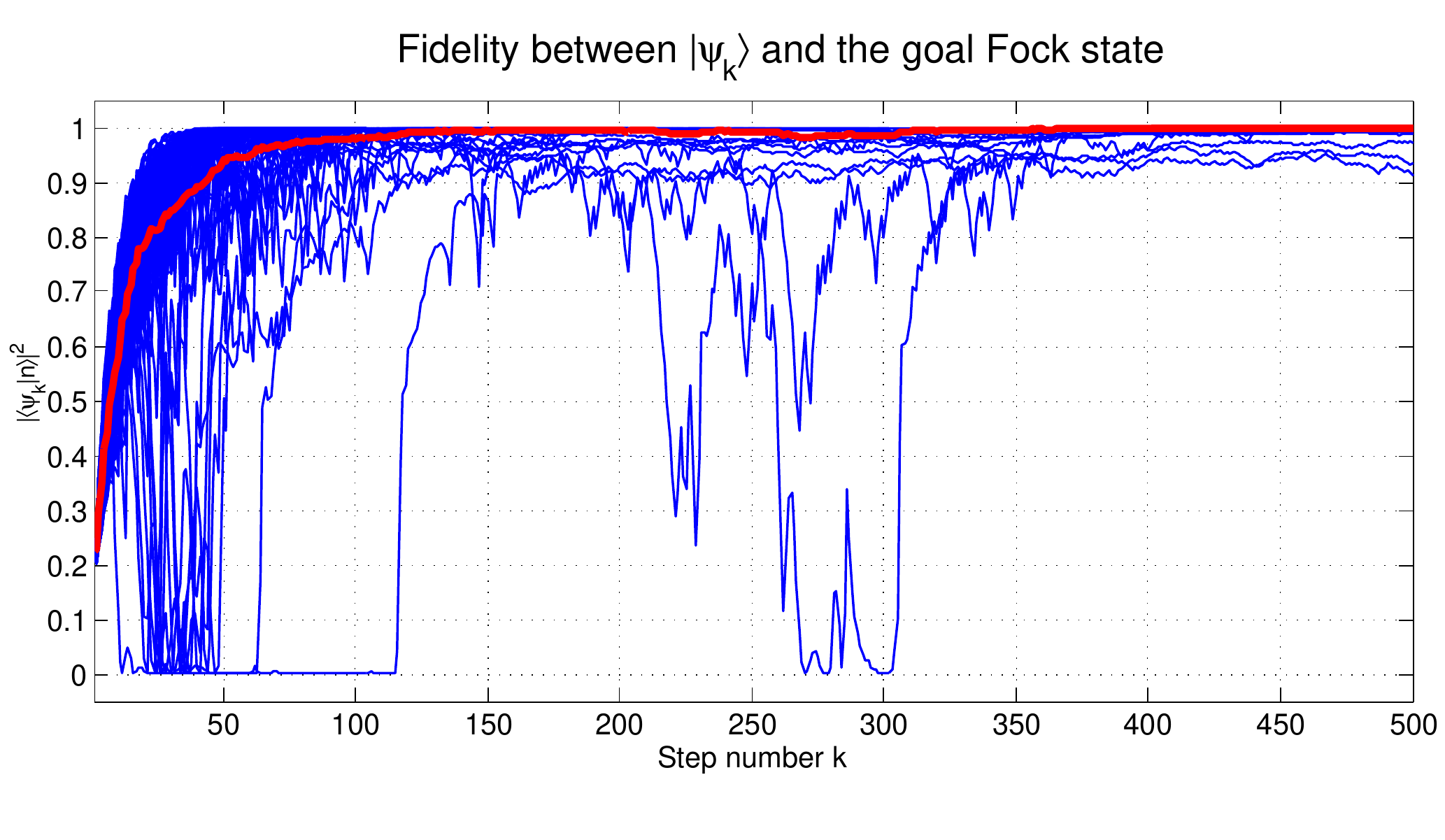}\\
  \caption{Simulation with a truncation to $20$ photons of the system and 9 photons of the filter for the feedback law~\eqref{eqn:control}; in blue$\left|\langle\bar n|\psi_k\rangle\right|^2 $ ($\bar n=3$) for each realization ; in red average over the 100 realizations of $\left|\langle\bar n|\psi_k\rangle\right|^2$.   }\label{fig:Traj}
  \includegraphics[trim = 0mm 85mm 0mm 90mm, clip, width=0.45\textwidth]{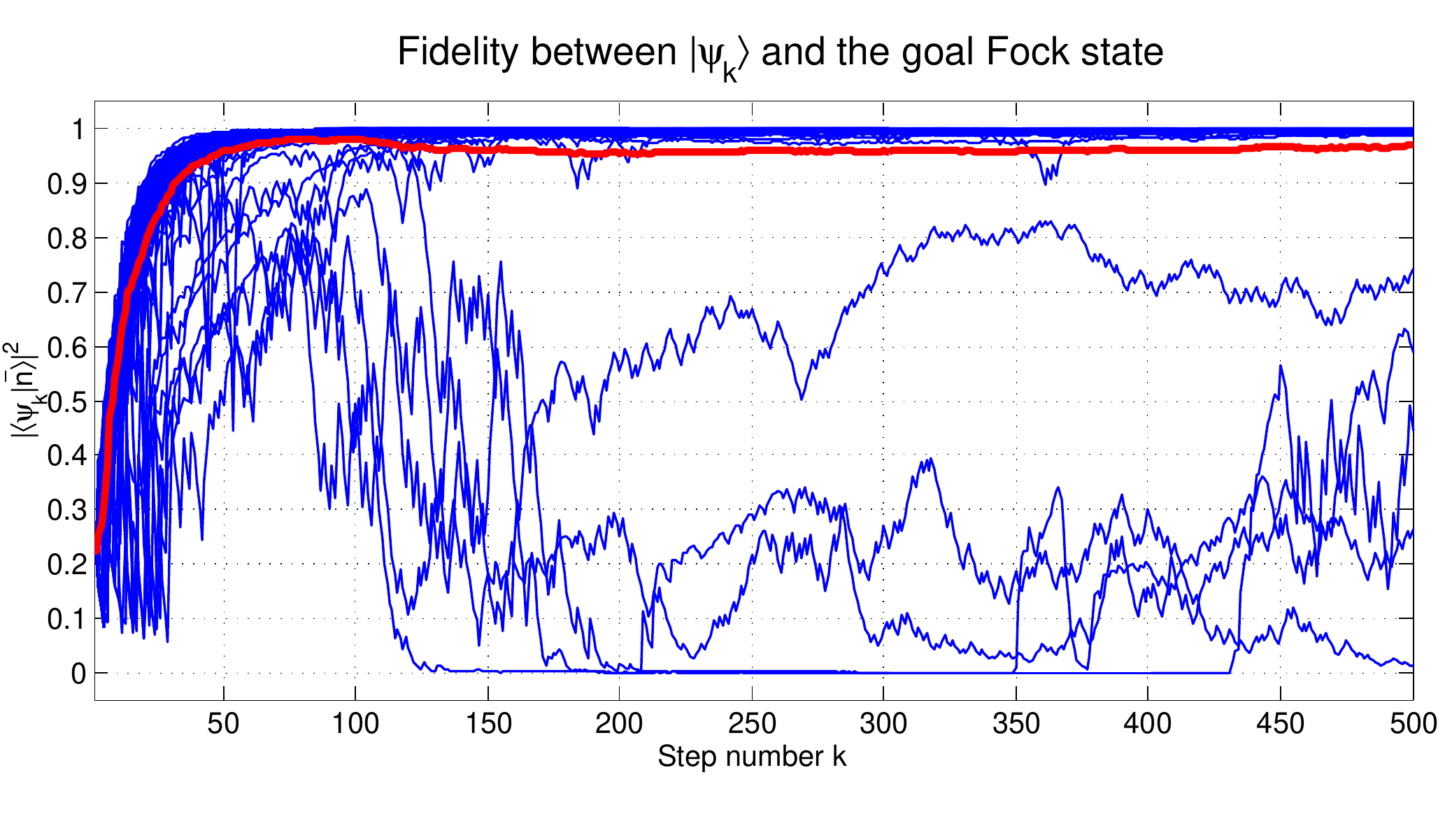}\\
 \caption{Simulation with a truncation to $20$ photons of the system and 9 photons of the filter for the "finite dimensional" feedback law~\eqref{eqn:controlPRA2009}; in blue$\left|\langle\bar n|\psi_k\rangle\right|^2 $ ($\bar n=3$) for each realization ; in red average over the 100 realizations of $\left|\langle\bar n|\psi_k\rangle\right|^2$.    }\label{fig:TrajPRA2009}
\end{figure}
\begin{figure}
\includegraphics[trim = 30mm 85mm 0mm 80mm, clip, width=0.6\textwidth]{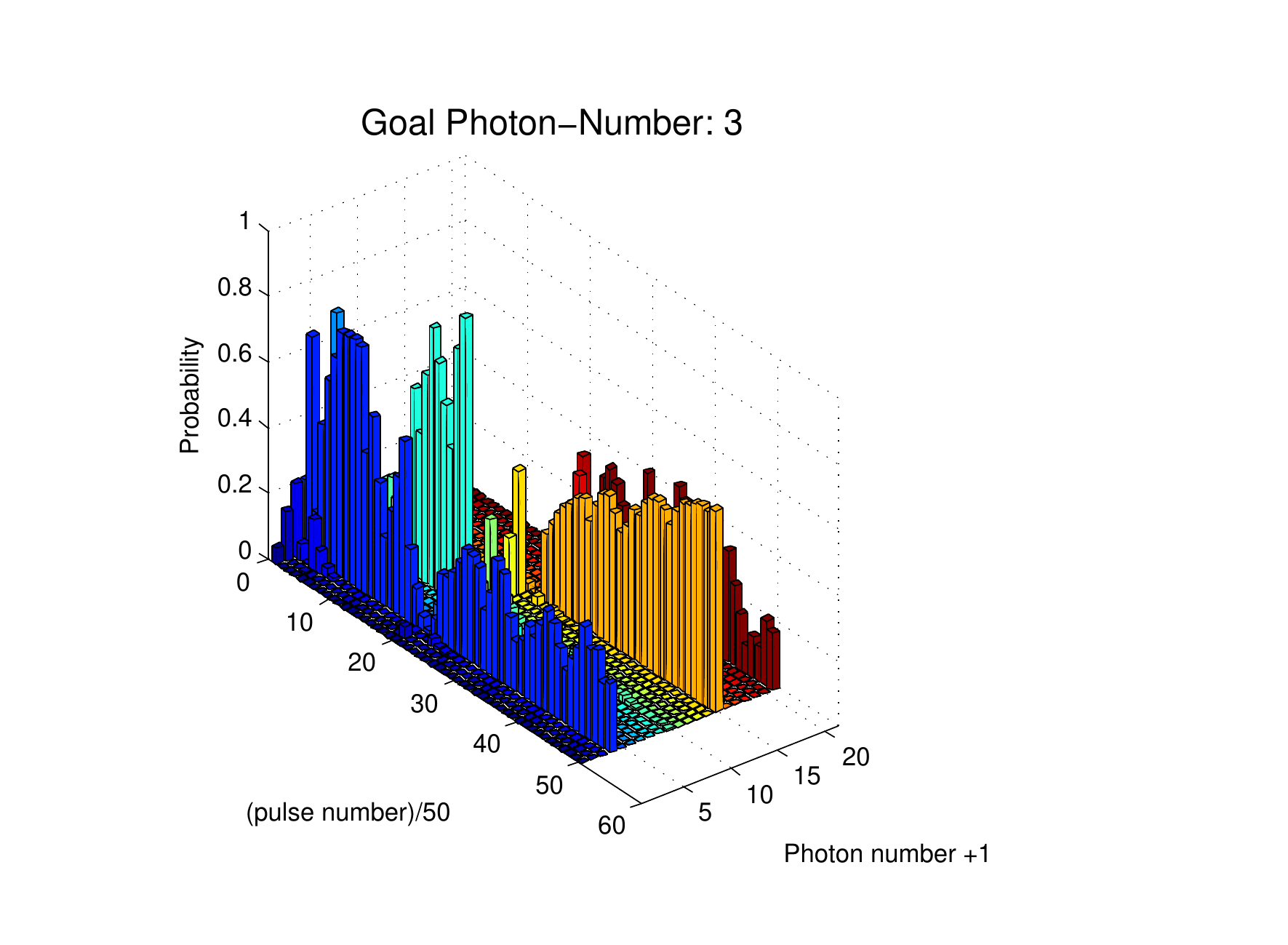}\\
\caption{An example of a trajectory of the finite-dimensional controller demonstrating escape to high photon numbers.}\label{fig:typTrajMassLoss}
\end{figure}
To illustrate Theorem~\ref{the:mainRes}, we performed closed-loop simulations of the controller designed using the finite-dimensional approximation~\cite{Mirrahimi2010} and the one in Theorem~\ref{the:mainRes}. Both simulations were performed on a system truncated to 21 photons. However the quantum filter (and therefore the controller) was truncated to 10 photons. 

The initial state was chosen to be the coherent state having  an average of $\bar n=3$ photons:  
$$
\ket{\psi_0}=e^{-\frac{\bar n}{2}}\sum_{n\geq 0} \sqrt{\tfrac{\bar n^n}{n!}} \ket{n}
$$
The measurement operators are $M_g=\cos\left(\sqrt{2}(N-\bar n)/5+\tfrac{\pi}{4}\right)$, $M_e=\sin\left(\sqrt{2}(N-\bar n)/5+\tfrac{\pi}{4}\right)$. We take $\bar\alpha=\tfrac{1}{10}$ and $\delta = (1/10(10+1))/2$ to ensure the Lyapunov function is strictly concave near the Fock states $\ket{m}$, $m\neq \bar{n}$. To compute the feedback law given by the minimisation~\eqref{eqn:control}, we approximate, for each step $k$, $[-\bar\alpha,+\bar\alpha]\ni\alpha_k\mapsto \bbE \left[V\left(\ket{\psi_{k+1}}\right)\big| \ket{\psi_k}\right]$ by the polynomial of degree two with the same first and second order derivatives at $\alpha_k=0$.  Figure~\ref{fig:Traj} shows  good convergence properties of such feedback strategy with an average asymptotic value of $\left|\langle\bar n|\psi\rangle\right|^2$  close to $1$. The remaining trajectories that do not converge to $\ket{\bar n}$ can be interpreted as the $\epsilon$ in theorem~\ref{the:mainRes}.

Figure~\ref{fig:TrajPRA2009} is devoted to similar simulations but with the  feedback law of~\cite{Mirrahimi2010,Dotsenko2009} based on a finite dimensional model:
\begin{equation}\label{eqn:controlPRA2009}
\alpha_k = 
\left\{
  \begin{array}{l}
     \bar\alpha \quad \hbox{if } \left|\langle\bar n|\psi_{k+1/2}\rangle\right|^2 \leq \tfrac{1}{10};
\\
   \tfrac{\left\langle \psi_{k+1/2} \big|\left[\kebr{\bar n}, a^\dag - a,\right]\psi_{k+1/2} \right\rangle}{4 \bar n +2} 
\quad \hbox{otherwise.}
  \end{array}
\right.
\end{equation}
The average asymptotic value of $\left|\langle\bar n|\psi\rangle\right|^2$  is then  around $0.95$ with this  "finite dimensional" feedback. Around $5\%$ of the trajectories do not converge towards $\ket{\bar n}$ and escape towards high photon numbers. Figure~\ref{fig:typTrajMassLoss} shows a typical example of such a trajectory which converges towards photon number 15 and 20.
\section{Conclusion}\label{sec:concl}
In this paper we examine the stabilization of a quantum optical cavity at a pre-specified photon number state $\ket{\bar{n}}$. In contrast with previous work, we designed a Lyapunov function on the entire infinite dimensional Hilbert space instead of using a truncation approximation. The Lyapunov function was chosen so that it is a strict Lyapunov function for the target state and the feedback consisted of a control that minimizes the expectation value of the Lyapunov function at each time-step. Simulations indicate that this feedback controller performs better than the one designed using the finite dimensional approximation.
\section{Acknowledgments}
The authors thank M. Brune, I.~Dotsenko, S.~Haroche and J.M. Raimond  for enlightening discussions and advices. 
\bibliographystyle{IEEEtran}
\bibliography{ref}
\end{document}